June 21, 2018

# CLOSE ENCOUNTERS WITH THE STIRLING NUMBERS OF THE SECOND KIND

Khristo N Boyadzhiev

Department of Mathematics and Statistics,
Ohio Northern University, Ada, Ohio 45810, USA
k-boyadzhiev@onu.edu

**Abstract** This is a short introduction to the theory of Stirling numbers of the second kind $S(m,k)$ from the point of view of analysis. It is written in the form of a historical survey. We tell the story of their birth in the book of James Stirling (1730) and show how they mature in the works of Johann Grünert (1843). We demonstrate their usefulness in several differentiation formulas. The reader can also see the connection of $S(m,k)$ to Bernoulli numbers, to Euler polynomials and to power sums.

## 1. Introduction

There are several sequences of special numbers which are very important in analysis, number theory, combinatorics, and even physics. Such are, for example, the Bernoulli numbers [1], [5], [9]. Equally important are also the Stirling numbers. We present here a short tutorial on the Stirling numbers of the second kind, including historical notes and some applications. We also show their relation to Bernoulli numbers, Euler polynomials, and sums of powers. The starting point is one interesting identity.
From time to time, in different books and articles we find the strange evaluation,

$$\sum_{k=0}^{n} \binom{n}{k}(-1)^k k^m = \begin{cases} 0 & \text{if } m<n \\ (-1)^n n! & \text{if } m=n \end{cases}. \qquad (1.1)$$

For instance, the following theorem (extending (1.1) a little bit) appeared in [12].

**Theorem**. *For any two real or complex numbers $x, y$ and for any positive integers $m$ and $n$,*

$$\sum_{k=0}^{n}\binom{n}{k}(-1)^k (xk+y)^m = \begin{cases} 0 & \text{if } m<n \\ (-1)^n x^n n! & \text{if } m=n \end{cases} \quad (1.2)$$

This, indeed, is a very strange result! Is this just a curious fact, or is there something bigger behind it? Also, what happens when $m > n$? The theorem obviously deserves further elaboration. Therefore, we want to fill this gap now and also provide some related historical information.

Identity (1.2) is not new. It appears in a more general form in the table of identities [7]. Namely, if $f(t) = c_0 + c_1 t + ... + c_m t^m$ is a polynomial of degree $m$, entry (Z.8) in [7] says that

$$\sum_{k=0}^{n}\binom{n}{k}(-1)^k f(k) = \begin{cases} 0 & \text{if } m<n \\ (-1)^n n! c_n & \text{if } m=n \end{cases} \quad (1.3)$$

which implies (1.2). As H. W. Gould writes on p. 82: "Relation (Z.8) is very useful; we have numerous interesting cases by choosing $f(t)$…". Identity (1.2) was later rediscovered by Ruiz [16], who proved it by induction.

Here is a simple observation (made also in [12]) – expanding the binomial $(xk + y)^m$ in (1.2) and changing the order of summation, we find that (1.2) is based on the more simple identity (1.1), (which, by the way, is entry (1.13) in [7]). In his paper [8] Henry Gould provides a nice and thorough discussion of identity (1.1), which he calls *Euler's formula*, as it appears in the works of Euler on n-th differences of powers. See also [17, pp. 18-19 and p. 48]. The mystery of identity (1.1) is revealed by its connection to the Stirling numbers. An old result in classical analysis (also discussed in [8]) says that

$$(-1)^n n! S(m,n) = \sum_{k=0}^{n}\binom{n}{k}(-1)^k k^m , \quad (1.4)$$

where $S(m,n)$ are the Stirling numbers of the second kind [2], [5], [8], [9], [11], [19]. They have the property $S(m,n) = 0$ when $m < n$, and $S(m,m) = 1$. We can define these numbers in combinatorial terms: $S(m,n)$ counts the number of ways to partition a set of $m$ elements into $n$ nonempty subsets. Thus we can read $S(m,n)$ as "m subset n". An excellent combinatorial treatment of the Stirling numbers can be found in [9]. Pippenger in the recent article [14] mentions, among other things, also their probabilistic interpretation. The alternative notation

$$S(m,n) = \begin{Bmatrix} m \\ n \end{Bmatrix}. \tag{1.5}$$

suggested in 1935 by the Serbian mathematician Jovan Karamata (see [9, p. 257]) fits very well with the combinatorial interpretation. A simple combinatorial argument ([9, p. 259]) provides the important recurrence

$$\begin{Bmatrix} m \\ n \end{Bmatrix} = n \begin{Bmatrix} m-1 \\ n \end{Bmatrix} + \begin{Bmatrix} m-1 \\ n-1 \end{Bmatrix}, \quad (n \leq m), \tag{1.6}$$

which parallels the well-known property of binomial coefficients [9, p158]

$$\binom{m}{n} = \binom{m-1}{n} + \binom{m-1}{n-1}. \tag{1.7}$$

For example, from (1.6) we can compute

$$\begin{Bmatrix} n+1 \\ n \end{Bmatrix} = n \begin{Bmatrix} n \\ n \end{Bmatrix} + \begin{Bmatrix} n \\ n-1 \end{Bmatrix} = n + \begin{Bmatrix} n \\ n-1 \end{Bmatrix} = n + (n-1) + \begin{Bmatrix} n-1 \\ n-2 \end{Bmatrix} \tag{1.8}$$

$$= n + (n-1) + (n-2) + ... + 1 = \frac{n(n+1)}{2},$$

i.e. $S(n+1, n) = \dfrac{n(n+1)}{2}$, and from (1.4) we find

$$\sum_{k=0}^{n} \binom{n}{k}(-1)^k k^{n+1} = \frac{(-1)^n n}{2}(n+1)!. \tag{1.9}$$

A proof of (1.6) using (1.7) and finite differences is presented in [8].

With the help of equation (1.4) we can fill the gap in (1.2) for $m > n$. Namely, we have

$$\sum_{k=0}^{n} \binom{n}{k}(-1)^k (xk+y)^m = (-1)^n n! \sum_{j=n}^{m} \binom{m}{j} x^j y^{m-j} S(j,n) \tag{1.10}$$

For the proof we just need to expand $(xk+y)^m$, change the order of summation and apply (1.4).

Identity (1.3) has a short and nice extension beyond polynomials. The following representation is true

$$\sum_{k=0}^{n}\binom{n}{k}(-1)^k f(k) = (-1)^n n! \sum_{m=0}^{\infty} c_m S(m,n) \tag{1.11}$$

for any $n \geq 0$ and any function $f(t) = c_0 + c_1 t + ...$ analytic on a disk with radius $R > n$.

Proof: We multiply (1.4) by $c_m$ and sum for $m$ from zero to infinity.

A combinatorial proof of (1.4) based on the combinatorial definition of $S(m,n)$ can be found in [5, pp. 204-205]; a proof based on finite differences in given in [11, p. 169] (see also [11, pp.177-178]). We shall present here two proofs of (1.1) and (1.4). For the first one we shall visit the birthplace of the Stirling numbers.

## 2. James Stirling and his Table

The name "Stirling numbers" comes from the Danish mathematician Niels Nielsen (1865-1931). On p. 68 of his book [13] Nielsen attributed these numbers to James Stirling, a Scottish mathematician (1692-1770), who worked on Newton series. At that time mathematicians realized the importance series expansion of functions and various techniques were gaining momentum.

Stirling studied in Oxford, then went to Italy for political reasons and became a professor of mathematics in Venice. In 1718 he published through Newton a paper titled "Methodus Differentialis Newtoniana Illustrata". In 1725 Stirling returned to England and in 1730 published his book "Methodus Differentialis" (The Method of Differences) [18]. The book was written in Latin, as most scientific books of that time were written. An annotated English translation was published recently by Ian Tweddle [19].

*Methodus Differentialis:*

SIVE

**TRACTATUS**

DE

**SUMMATIONE**

ET

**INTERPOLATIONE**

SERIERUM INFINITARUM.

AUCTORE *JACOBO STIRLING*, R.S.S.

Figure 1. Part of the front page of Stirling's book

At the beginning of his book Stirling studied carefully the coefficients $A_m^n$ in the representations

$$z^m = A_1^m z + A_2^m z(z-1) + A_3^m z(z-1)(z-2) + ... + A_m^m z(z-1)...(z-m+1) \ , \qquad (2.1)$$

where $m = 1, 2, ...$ . On p. 8 he presented a table containing many of these coefficients

<div style="text-align:center">*Tabulam priorem.*</div>

| 1 | 1 | 1 | 1 | 1 | 1 | 1 | 1 | 1 | &c. |
|---|---|---|---|---|---|---|---|---|---|
|   | 1 | 3 | 7 | 15 | 31 | 63 | 127 | 255 | &c. |
|   |   | 1 | 6 | 25 | 90 | 301 | 966 | 3025 | &c. |
|   |   |   | 1 | 10 | 65 | 350 | 1701 | 7770 | &c. |
|   |   |   |   | 1 | 15 | 140 | 1050 | 6951 | &c. |
|   |   |   |   |   | 1 | 21 | 266 | 2646 | &c. |
|   |   |   |   |   |   | 1 | 28 | 461 | &c. |
|   |   |   |   |   |   |   | 1 | 36 | &c. |
|   |   |   |   |   |   |   |   | 1 | &c. |
|   |   |   |   |   |   |   |   |   | &c. |

Figure 2. Stirling's First Table

In the table $m$ changes horizontally, left to right, and $n$ changes vertically, from top to bottom. Therefore, by following the columns of the table we find

$$z = z, \qquad (2.2)$$
$$z^2 = z + z(z-1),$$
$$z^3 = z + 3z(z-1) + z(z-1)(z-2),$$
$$z^4 = z + 7z(z-1) + 6z(z-1)(z-2) + z(z-1)(z-2)(z-3),$$
$$z^5 = z + 15z(z-1) + 25z(z-1)(z-2) + 10z(z-1)(z-2)(z-3) + z(z-1)(z-2)(z-3)(z-4),$$
etc.

The coefficients $A_k^m$ are exactly the numbers which we call today Stirling numbers of the second kind. For completeness, we add to this sequence also $A_0^0 = 1$, $A_0^m = 0$ when $m > 0$. The following is true.

**Theorem A**. *Let the coefficients $A_n^m$ be defined by the expansion (2.1). Then*

$$A_n^m = \frac{1}{n!} \sum_{k=0}^{n} \binom{n}{k} (-1)^{n-k} k^m , \qquad (2.3)$$

*and the right hand side is zero, when $n > m$.*

It is not obvious that (2.1) implies (2.3). For the proof of the theorem we need some preparation.

Stirling's technique for computing this table is presented on pp. 24-29 in Ian Tweddle's translation [19]. As Ian Tweddle comments on p. 171 in [19], had Stirling known the relation (1.6), the computation of the table would have been much easier.

## 3. Newton series and finite differences

In the seventeenth and eighteenth centuries it was considered important to expand functions in Newton series. That is,

$$f(z) = \sum_{k=0}^{\infty} a_k z(z-1)(z-2)\ldots(z-k+1) \qquad (3.1)$$

$$= a_0 + a_1 z + a_2 z(z-1) + a_3 z(z-1)(z-2) + \ldots .$$

The theory of Newton series is similar to the theory of Taylor series. First of all, one needs to find a formula for the coefficients $a_k$. In the case of Taylor series the function $f(z)$ is expanded on the power polynomials $p_k(z) = z^k$, $k = 0,1,\ldots$, and the coefficients are expressed in terms of the higher derivatives of the function evaluated at zero. In the case of Newton series, instead of derivatives one needs to use finite differences. This is suggested by the very form of the series, as the function is expanded now on the difference polynomials $P_0(z) = 1, P_1(z) = z,\ldots, P_k(z) = z(z-1)\ldots(z-k+1),\ldots$, which constitute another basis in the space of polynomials.

For a given function $f(z)$ we set

$$\Delta f(z) = f(z+1) - f(z) . \qquad (3.2)$$

Then

$$\Delta^2 f(z) = \Delta(\Delta f)(z) = f(z+2) - 2f(z+1) + f(z), \tag{3.3}$$

$$\Delta^3 f(z) = \Delta(\Delta^2 f)(z) = f(z+3) - 3f(z+2) + 3f(z+1) - f(z),$$

etc.

We notice the binomial coefficients appearing here with alternating signs. Following this pattern we arrive at the representation

$$\Delta^n f(z) = \sum_{k=0}^{n} \binom{n}{k}(-1)^{n-k} f(z+k). \tag{3.4}$$

In particular, with $z = 0$,

$$\Delta^n f(0) = \sum_{k=0}^{n} \binom{n}{k}(-1)^{n-k} f(k). \tag{3.5}$$

This formula will be used now to compute the coefficients $a_k$ in the Newton series (3.1). With $z = 0$ we see that $a_0 = f(0)$. A simple computation shows that $\Delta P_k(z) = k P_{k-1}(z), k \geq 1$ and so

$$\Delta f(z) = a_1 + 2a_2 z + 3a_3 z(z-1) + \ldots , \tag{3.6}$$

which yields $a_1 = \Delta f(0)$. Also,

$$\Delta^2 f(z) = 2a_2 + 2 \cdot 3 a_3 z + 3 \cdot 4 z(z-1) + \ldots \tag{3.7}$$

and $2a_2 = \Delta^2 f(0)$. Continuing this way we find $3!a_3 = \Delta^3 f(0), 4!a_4 = \Delta^4 f(0)$, etc. The general formula is $k!a_k = \Delta^k f(0)$, $k = 0,1,\ldots$. Thus (3.1) becomes

$$f(z) = \sum_{k=0}^{\infty} \frac{\Delta^k f(0)}{k!} z(z-1)(z-2)\ldots(z-k+1). \tag{3.8}$$

*Proof of Theorem A.* Take $f(z) = z^m$ in (3.8) to obtain (in view of (3.5))

$$z^m = \sum_{n=0}^{\infty} \left\{ \frac{1}{n!} \sum_{k=0}^{n} \binom{n}{k}(-1)^{n-k} k^m \right\} z(z-1)(z-2)\ldots(z-n+1) . \tag{3.9}$$

Comparing this to (2.1) yields the representation (2.3). Note also that the series in (3.9) truncates, as on the left-hand side we have a polynomial of degree $m$. The summation at the right-hand side stops with $n = m$, and

$$\frac{1}{m!} \sum_{k=0}^{m} \binom{m}{k} (-1)^{m-k} k^m = 1 \ . \tag{3.10}$$

This is not the end of the story, however. In the curriculum vitae of the Stirling numbers there is one other most remarkable event.

## 4. Grünert's polynomials

Amazingly, the same Stirling numbers appeared again, one hundred years later, in a very different setting. They appeared in the work [10] of the German mathematician Johann August Grünert (1797-1872), professor at the University of Greifswald, Germany. He taught there from 1833 until his death. Grünert, a student of Pfaff and Gauss, was interested in many topics, not only in mathematics, but also in physics. He wrote a number of books on such diverse subjects as conic sections, the loxodrome, optics, and the solar eclipse. Some of his books, including Optische Untersuchungen (Studies in Optics) and Theorie der Sonnenfinsternisse (Theory of the Solar Eclipse) are available now as Google books on the Internet. In 1841 Grünert started to edit and publish the highly respected Archiv der Mathematik und Physik (known also as "Grünert's Archiv"). His biography, written by his student Maximus Curtze, appeared in volume 55 of that journal.

Grünert came to the numbers $S(m,n)$ by repeatedly applying the operator $x\frac{d}{dx}$ to the exponential function $e^x$. This procedure generates a sequence of polynomials

$$x\frac{d}{dx}e^x = xe^x,$$

$$\left(x\frac{d}{dx}\right)^2 e^x = (x^2 + x)e^x,$$

$$\left(x\frac{d}{dx}\right)^3 e^x = (x^3 + 3x^2 + x)e^x,$$

...

$$\left(x\frac{d}{dx}\right)^m e^x = (B_0^m + B_1^m x + B_2^m x^2 + ... + B_m^m x^m)e^x, \tag{4.1}$$

etc, with certain coefficients $B_k^m$. We shall see that these coefficients are exactly the Stirling numbers of the second kind. This fact follows from the theorem.

**Theorem B** (Grünert). *Let the coefficients $B_k^m$ be defined by equation (4.1). Then*

$$B_n^m = \frac{1}{n!}\sum_{k=0}^{n}\binom{n}{k}(-1)^{n-k} k^m . \tag{4.2}$$

**Conclusion.** From (4.2), (2.3) and (1.4) we conclude that $A_n^m = B_n^m = S(m,n)$.

*Proof of Theorem B.* From the expansion

$$e^x = \sum_{k=0}^{\infty}\frac{x^k}{k!} , \tag{4.3}$$

we find

$$\left(x\frac{d}{dx}\right)^m e^x = \sum_{k=0}^{\infty}\frac{k^m x^k}{k!} , \tag{4.4}$$

for $m = 0,1,...$ , and from (4.1),

$$B_0^m + B_1^m x + B_2^m x^2 + ... + B_m^m x^m = e^{-x}\sum_{k=0}^{\infty}\frac{k^m x^k}{k!} = \left\{\sum_{j=0}^{\infty}\frac{(-1)^j x^j}{j!}\right\}\left\{\sum_{k=0}^{\infty}\frac{k^m x^k}{k!}\right\}. \tag{4.5}$$

Multiplying the two power series on the right-hand side yields

$$B_0^m + B_1^m x + B_2^m x^2 + ... + B_m^m x^m = \sum_{n=0}^{\infty} x^n \left\{\frac{1}{n!}\sum_{k=0}^{n}\binom{n}{k}(-1)^{n-k} k^m\right\}, \tag{4.6}$$

and again, comparing coefficients we conclude that the series on the right-hand side is finite and (4.2) holds. The theorem is proved!

## 5. Intermediate summary and the exponential polynomials

So far we have established two important mathematical facts. Namely,

**Fact 1** The coefficients $A_n^m$ defined by the representation (2.1) are the same as the coefficients $B_n^m$ defined by equation (4.1) and

$$A_n^m = B_n^m = \frac{1}{n!}\sum_{k=0}^{n}\binom{n}{k}(-1)^{n-k} k^m . \tag{5.1}$$

These are the Stirling numbers of the second kind with present day notation in analysis $S(m,n)$ ([5], [19]). At that $S(m,n) = 0$ when $m < n$, and $S(n,n) = 1$. In particular, this proves (1.1).

**Fact 2** The Stirling numbers of the second kind are used in combinatorics, often with the notation $S(m,n) = \left\{{m \atop n}\right\}$ (see [2], [3], [4], [9], [14]). They represent the number of ways by which a set of $m$ elements can be partitioned into $n$ non-empty subsets. Thus $\left\{{m \atop n}\right\}$ is naturally defined for $n \leq m$ and $\left\{{n \atop n}\right\} = 1$. When $m < n$, $\left\{{m \atop n}\right\} = 0$. The numbers $\left\{{m \atop n}\right\}$ equal $A_n^m$ because they satisfy (1.4) as proven in [5].

The polynomials

$$\varphi_n(x) = S(n,0) + S(n,1)x + \ldots + S(n,n)x^n, \tag{5.2}$$

$n = 0, 1, \ldots$, appearing in Grünert's work are called exponential polynomials. They were later rediscovered and used by several authors. These polynomials are defined by equation (4.1), i.e.

$$\varphi_n(x) = e^{-x}\left(x\frac{d}{dx}\right)^m e^x \tag{5.3}$$

or, by the generating function (see [2]),

$$e^{x(e^t-1)} = \sum_{n=0}^{\infty} \frac{\varphi_n(x)}{n!} t^n. \tag{5.4}$$

Here are the first five of them

$$\begin{aligned}
\varphi_0(x) &= 1 \\
\varphi_1(x) &= x \\
\varphi_2(x) &= x^2 + x \\
\varphi_3(x) &= x^3 + 3x^2 + x \\
\varphi_4(x) &= x^4 + 6x^3 + 7x^2 + x
\end{aligned} \tag{5.5}$$

A short review of these polynomials is given in [2]. Interesting notes on the Stirling numbers can be found also in Gould's paper [8]. As Gould mentions in [8], these numbers appeared in the works of Euler, in the form $S(m,n)n!$. We shall discuss this in the next section.

Replacing $x$ by $ax$ in (5.3), where $a$ is any constant, we see that (5.3) can be written as

$$\left(x\frac{d}{dx}\right)^m e^{ax} = \varphi_n(ax)e^{ax}. \tag{5.6}$$

This form is useful in some computations.

The Stirling numbers of the second kind are closely related to the Stirling numbers of the first kind, i.e. $s(m,n)$, see [2], [5], [9], [19]. These numbers also appeared in Stirling's book [18]. For completeness we shall mention them in the last section.

## 6. The exponential generating function for $S(m,n)$

For any integer $n \geq 0$ let us expand the function $f(x) = (e^x - 1)^n$ in Maclaurin series, i.e.

$$f(x) = \sum_{m=0}^{\infty} \frac{f^{(m)}(0)}{m!} x^m. \tag{6.1}$$

For this purpose we first write

$$(e^x - 1)^n = \sum_{k=0}^{n} \binom{n}{k}(-1)^{n-k} e^{kx}, \tag{6.2}$$

and then, according to (5.1) we compute $f^{(m)}(0)$,

$$\left(\frac{d}{dx}\right)^m (e^x - 1)^n \bigg|_{x=0} = \sum_{k=0}^{n} \binom{n}{k}(-1)^{n-k} k^m = n! S(m,n).$$

Therefore,

$$\frac{1}{n!}(e^x - 1)^n = \sum_{m=0}^{\infty} S(m,n) \frac{x^m}{m!}. \tag{6.3}$$

This is the exponential generating function for the Stirling numbers of the second kind $S(m,n)$. The summation, in fact, goes for $m \geq n$, as $S(m,n) = 0$ when $m < n$. Equation (6.3) is often used for the initial definition of $S(m,n)$.

## 7. Euler and the derivatives game

Let $|x| < 1$. We want to show that the numbers $S(m,n)$ naturally appear in the derivatives

$$\left(x\frac{d}{dx}\right)^m \frac{1}{1-x} = \sum_{n=0}^{\infty} n^m x^n, \tag{7.1}$$

where $m = 0, 1, \ldots$. To show this we first write

$$\frac{1}{1-x} = \int_0^{\infty} e^{-(1-x)t} dt = \int_0^{\infty} e^{xt} e^{-t} dt. \tag{7.2}$$

Then, in view of (5.6),

$$\left(x\frac{d}{dx}\right)^m \frac{1}{1-x} = \int_0^{\infty} \varphi_m(xt) e^{xt} e^{-t} dt = \sum_{n=0}^{m} S(m,n) x^n \int_0^{\infty} t^n e^{-(1-x)t} dt \tag{7.3}$$

$$= \sum_{n=0}^{m} S(m,n) x^n \frac{n!}{(1-x)^{n+1}} = \frac{1}{1-x} \sum_{n=0}^{m} S(m,n) n! \left(\frac{x}{1-x}\right)^n.$$

For the third equality we use the well-known formula (Laplace transform of $t^{\alpha}$)

$$\frac{\Gamma(\alpha+1)}{s^{\alpha+1}} = \int_0^{\infty} t^{\alpha} e^{-st} dt. \tag{7.4}$$

Introducing the polynomials

$$\omega_m(z) = \sum_{n=0}^{m} S(m,n) n! z^n, \tag{7.5}$$

we can write (7.3) in the form

$$\sum_{n=0}^{\infty} n^m x^n = \frac{1}{1-x} \omega_m\left(\frac{x}{1-x}\right). \tag{7.6}$$

Thus we have

$$\omega_0(x) = 1$$
$$\omega_1(x) = x$$
$$\omega_2(x) = 2x^2 + x \qquad (7.7)$$
$$\omega_3(x) = 6x^3 + 6x^2 + x$$
$$\omega_4(x) = 24x^4 + 36x^3 + 14x^2 + x ,$$

etc. The polynomials $\omega_n$ can be seen on p. 389, Chapter VII, in Euler's book [6].

Figure 3. The geometric polynomials in Euler's work

Essentially, Euler obtained these polynomials by computing the derivatives (7.1) directly. We shall see now how all this can be done in terms of exponentials.

Here is a good exercise. Let us expand the function

$$f(t) = \frac{1}{\mu e^{\lambda t} + 1} \qquad (7.8)$$

in Maclaurin series ($\lambda, \mu$ are two parameters). We need to find the higher derivatives of $f$ at zero. Assuming for the moment that $|\mu e^{\lambda t}| < 1$ we use the expansion

$$\frac{1}{\mu e^{\lambda t}+1}=\frac{1}{1-(-\mu e^{\lambda t})}=\sum_{n=0}^{\infty}(-\mu)^{n}e^{\lambda t n}. \tag{7.9}$$

From this

$$\left(\frac{d}{dt}\right)^{m}\frac{1}{\mu e^{\lambda t}+1}=\lambda^{m}\sum_{n=0}^{\infty}(-\mu)^{n}n^{m}e^{\lambda t n} \tag{7.10}$$

and in view of (7.6)

$$\left(\frac{d}{dt}\right)^{m}\frac{1}{\mu e^{\lambda t}+1}=\frac{\lambda^{m}}{\mu e^{\lambda t}+1}\omega_{m}\left(\frac{-\mu e^{\lambda t}}{\mu e^{\lambda t}+1}\right), \tag{7.11}$$

so that

$$\left(\frac{d}{dt}\right)^{m}\frac{1}{\mu e^{\lambda t}+1}\bigg|_{t=0}=\frac{\lambda^{m}}{\mu+1}\omega_{m}\left(\frac{-\mu}{\mu+1}\right), \tag{7.12}$$

which yields the desired representation

$$\frac{1}{\mu e^{\lambda t}+1}=\frac{1}{\mu+1}\sum_{m=0}^{\infty}\lambda^{m}\omega_{m}\left(\frac{-\mu}{\mu+1}\right)\frac{t^{m}}{m!}. \tag{7.13}$$

In particular, with $\lambda=\mu=1$,

$$\frac{2}{e^{t}+1}=\sum_{m=0}^{\infty}\omega_{m}\left(\frac{-1}{2}\right)\frac{t^{m}}{m!}=\sum_{m=0}^{\infty}\left\{\sum_{n=0}^{m}S(m,n)n!\frac{(-1)^{n}}{2^{n}}\right\}\frac{t^{m}}{m!}. \tag{7.14}$$

The polynomials $\omega_{m}$ appeared in the works of Euler, but they do not carry his name. In [8] and [17], they are used to evaluate the series on the right-hand side of (6.1) in terms of Stirling numbers. These polynomials were studied in [4] and called *geometric polynomials,* because of their obvious relation to the geometric series. It was shown in [4] that $\omega_{m}$ participate in a certain *series transformation formula*. In [3] the geometric polynomials were used to compute the derivative polynomials for $\tan x$ and $\sec x$.

One can write (7.6) in the form

$$\sum_{n=0}^{\infty}n^{m}x^{n}=\frac{A_{m}(x)}{(1-x)^{m+1}} \tag{7.15}$$

where $A_m$ are polynomials of degree $m$. These polynomials are known today as Eulerian polynomials and their coefficients are the Eulerian numbers [5], [9].

At the same time, there is a sequence of interesting and important polynomials carrying the name Euler polynomials. These are the polynomials $E_m(x)$, $m = 0,1,...$, defined by the generating function

$$\frac{2e^{xt}}{e^t+1} = \sum_{m=0}^{\infty} E_m(x)\frac{t^m}{m!} .\tag{7.16}$$

Using (7.14) we write (as in (4.5))

$$\frac{2e^{xt}}{e^t+1} = \left\{\sum_{n=0}^{\infty} \frac{x^n t^n}{n!}\right\}\left\{\sum_{k=0}^{\infty} \omega_k\left(\frac{-1}{2}\right)\frac{t^k}{k!}\right\} = \sum_{m=0}^{\infty} \frac{t^m}{m!}\left\{\sum_{k=0}^{m}\binom{m}{k}\omega_k\left(\frac{-1}{2}\right)x^{m-k}\right\} .\tag{7.17}$$

Comparing (7.16) and (7.17) yields

$$E_m(x) = \sum_{k=0}^{m}\binom{m}{k}\omega_k\left(\frac{-1}{2}\right)x^{m-k} ,\tag{7.18}$$

with

$$E_m(0) = \omega_m\left(\frac{-1}{2}\right) = \sum_{n=0}^{m} S(m,n)\, n!\,\frac{(-1)^n}{2^n} .\tag{7.19}$$

## 8. Relation to Bernoulli numbers

The Bernoulli numbers $B_m$, $m = 0,1,...$, can be defined by the generating function

$$\frac{t}{e^t-1} = \sum_{m=0}^{\infty} B_m \frac{t^m}{m!} , \quad |t| < 2\pi \tag{8.1}$$

([1], [5], [9]). From this

$$B_m = \left(\frac{d}{dt}\right)^m \left.\frac{t}{e^t-1}\right|_{t=0} .\tag{8.2}$$

It is tempting to evaluate these derivatives at zero by using the Leibnitz rule for the product $t \cdot \frac{1}{e^t-1}$ and formula (7.12) with $\mu = -1$, $\lambda = 1$. This will not work, though, because the denominator $\mu+1$ on the right-hand side becomes zero. To find a relation between the Bernoulli and Stirling numbers we shall use a simple trick and the generating function (6.3). Writing $t = \ln e^t = \ln(1+(e^t-1))$ we have for $t$ small enough

$$\frac{t}{e^t-1} = \frac{\ln(1+(e^t-1))}{e^t-1} = \sum_{n=0}^{\infty} \frac{(-1)^n}{n+1}(e^t-1)^n \qquad (8.3)$$

$$= \sum_{n=0}^{\infty} \frac{(-1)^n}{n+1}\left\{ n! \sum_{m=n}^{\infty} S(m,n) \frac{t^m}{m!} \right\} = \sum_{m=0}^{\infty} \frac{t^m}{m!}\left\{ \sum_{n=0}^{m} (-1)^n \frac{n!}{n+1} S(m,n) \right\}.$$

Comparing this to (8.1) we find for $m = 0,1,...$,

$$B_m = \sum_{n=0}^{m} (-1)^n \frac{n!}{n+1} S(m,n) . \qquad (8.4)$$

## 9. Sums of powers

The Bernoulli numbers historically appeared in the works of the Swiss mathematician Jacob Bernoulli (1654-1705) who evaluated sums of powers of consecutive integers (see [1], [9])

$$1^m + 2^m + ... + (n-1)^m = \frac{1}{m+1} \sum_{k=0}^{m} \binom{m+1}{k} B_k n^{m+1-k}, \qquad (9.1)$$

for any $m \geq 0, n \geq 1$. This is the famous Bernoulli formula. It is interesting to see that sums of powers can also be evaluated directly in terms of Stirling numbers of the second kind. In order to do this, we invert the representation (1.4). i.e.

$$n^m = \sum_{k=0}^{n} \binom{n}{k} S(m,k) k!. \qquad (9.2)$$

This inversion is a property of the binomial transform [15]. Given a sequence $\{a_k\}$, its *binomial transform* $\{b_k\}$ is the sequence defined by

$$b_n = \sum_{k=0}^{n} \binom{n}{k} a_k, \qquad (9.3)$$

and the inversion formula is

$$a_n = \sum_{k=0}^{n} \binom{n}{k} (-1)^{n-k} b_k . \qquad (9.4)$$

Next, from (9.2),

$$1^m + 2^m + ... + n^m = \sum_{p=1}^{n} \left\{ \sum_{k=0}^{p} \binom{p}{k} S(m,k) k! \right\} \tag{9.5}$$

$$= \sum_{k=0}^{n} S(m,k) k! \left\{ \sum_{p=k}^{n} \binom{p}{k} \right\},$$

by changing the order of summation. Using now the well known identity

$$\sum_{p=k}^{n} \binom{p}{k} = \binom{n+1}{k+1}, \tag{9.6}$$

we finally obtain

$$1^m + 2^m + ... + n^m = \sum_{k=0}^{n} \binom{n+1}{k+1} S(m,k) k! , \tag{9.7}$$

which is the desired representation.

## 10. Stirling numbers of the first kind

Inverting equation (2.1) one has

$$z(z-1)...(z-m+1) = \sum_{k=0}^{m} s(m,k) z^k \tag{10.1}$$

where the coefficients $s(m,k)$ are called Stirling numbers of the first kind. The following inversion property is true

$$\sum_{k=0}^{\infty} S(m,k) s(k,n) = \delta_{mn} = \begin{cases} 0 & m \neq n \\ 1 & m = n \end{cases}. \tag{10.2}$$

Stirling's book [18], [19] contains a second table of coefficients

**Figure 4. Stirling's Second table**

The coefficients here come from the representation ($m = 1, 2, ...$)

$$\frac{1}{z^{m+1}} = \sum_{k=0}^{\infty} \frac{\sigma(m+k, m)}{z(z+1)...(z+m+k)}, \qquad (10.3)$$

following the columns of the table. The numbers $\sigma(m,k)$ are called today the Stirling cycle numbers or the unsigned Stirling numbers of the first kind [5], [9]. An often used notation is

$$\sigma(m,k) = \begin{bmatrix} m \\ k \end{bmatrix}. \qquad (10.4)$$

One has $s(m,k) = (-1)^{m-k} \sigma(m,k)$.

More properties, combinatorial interpretation, details and generating functions can be found in the excellent books [5], [9], [11].